\documentclass{amsart}

\usepackage{ stmaryrd }
\usepackage{amssymb}
\usepackage{amsthm}
\usepackage{amsmath}
\usepackage{amsmath,amsfonts,epsfig}
\usepackage{mathptmx}

\usepackage{amsmath,amsfonts,latexsym,amssymb}
\usepackage{mathrsfs}
\usepackage{verbatim}
\usepackage[utf8]{inputenc}
\usepackage[T1]{fontenc}
\usepackage{ae,aecompl}
\usepackage{braket}
\usepackage{color}
\usepackage{euscript}

\renewcommand{\bold}[1]{\medskip \noindent {\bf \boldmath #1
                        }\nopagebreak[4]}

\newcommand{\id}{\operatorname{id}}

\renewcommand{\Im}{\operatorname{Im}}

\newcommand{\sech}{\operatorname{sech}}

\newcommand{\Teich}{\operatorname{Teich}}

\newtheorem{theorem}{Theorem}[section]

\newtheorem{lemma}[theorem]{Lemma}

\newtheorem{corollary}[theorem]{Corollary}

\newcommand{\dome}{\operatorname{Dome}}
\renewcommand{\hbar}{\bar{{\mathbb H}}^3}

\newcommand{\CC}{\mathbb C}

\newcommand{\R}{\mathbb R}

\newcommand{\Z}{\mathbb Z}
\newcommand{\Half}{{\mathbb H}}
\newcommand{\Hp}{{\mathbb H}^2}
\newcommand{\Hs}{{\mathbb H}^3}

\newcommand{\psl}{{\rm PSL}(2,\CC)}
\newcommand{\pslr}{{\rm PSL}(2,\R)}

\newcommand{\CP}{{\operatorname{\mathbb{CP}}}}

\newcommand{\cotan}{\operatorname{cotan}}

\def\eproof{$\Box$ \medskip}

\renewcommand\marginpar[1]{} 

\begin{document}

\title[Schwarzian Bounds on Bending in Hyperbolic 3-Manifolds]{Schwarzian Bounds on Bending in Hyperbolic 3-Manifolds}
\author[Bridgeman]{Martin Bridgeman}
\address{Martin Bridgeman, Department of Mathematics, Boston College
Chestnut Hill, Ma 02467}
\thanks{Bridgeman's research is supported by NSF grant DMS-2405291 and Simons Foundation Travel Grant MPS-TSM-00007295}
\author[Tee]{Ming Hong Tee}
\address{Ming Hong Tee, Department of Mathematics, Boston College,
Chestnut Hill, Ma 02467}
\date{\today}

\maketitle
\begin{abstract} The Schwarzian derivative provides a classical analytic measure of how far a holomorphic map of the disk is from being M\"obius, with Nehari's bounds giving sharp criteria for univalence. Independently, Thurston introduced a geometric parametrization of locally univalent maps via bending measured laminations on the hyperbolic plane, capturing deviation from roundness in hyperbolic three-space. While both approaches quantify the same phenomenon, their precise relationship has remained only implicit.
In this paper we establish explicit quantitative bounds relating the Schwarzian norm $\|Sf\|_\infty$ and the bending norm $\|\beta_f\|_L$. In particular, for univalent maps with $\|Sf\|_\infty< 1/2$, we show that $\|\beta_f\|_L$ is controlled by an elementary function $B_L(\|Sf\|_L)$ that we compute explicitly. As an application, we obtain new effective bounds on the bending laminations of quasifuchsian manifolds in terms of the Teichm\"uller distance between their conformal boundary components. Our results sharpen the analytic-geometric correspondence between the Schwarzian derivative and hyperbolic geometry showing that just as small Schwarzian norm forces injectivity, it also forces controlled bending of convex hull boundaries.\end{abstract}
\section{Introduction}

The Schwarzian derivative has long played a central role in the study of univalent and locally univalent maps. For a holomorphic map 
\( f:\Delta \to \widehat{\mathbb{C}} \), the condition that the Schwarzian derivative \( Sf \) be small provides a powerful analytic criterion for injectivity. 
Classical results of Nehari (see \cite{nehari}) show that bounds on \( \|Sf\|_\infty \) determine when \( f \) is univalent: if \( \|Sf\|_\infty < \tfrac{1}{2} \), then 
\( f \) is univalent, while every univalent map satisfies \( \|Sf\|_\infty < \tfrac{3}{2} \). These estimates are fundamental to the analytic approach 
to Teichm\"uller theory, underlying Bers's model of Teichm\"uller space as a bounded domain in the space of quadratic differentials.  

On the other hand, Thurston introduced a geometric parametrization of locally univalent maps by associating to \( f \) a bending measured 
lamination \( \beta_f \) on the hyperbolic plane. This description connects the analysis of the Schwarzian to the geometry of convex hulls in hyperbolic 
three-space. The bending lamination provides a natural quantitative measure of how far the image of \( f \) is from being round, and has proved 
central in understanding the geometry of quasifuchsian groups and convex cores of hyperbolic 3-manifolds.  

Both perspectives, Nehari's analytic criterion and Thurston's geometric parametrization, give complementary ways to measure deviation from 
M\"obius maps. It is therefore natural to ask: 
\emph{how are the analytic and geometric measures related?} Specifically, how can one control the bending norm \( \|\beta_f\|_L \) in terms of the 
Schwarzian norm \( \|Sf\|_\infty \)?  

Implicit compactness arguments (Bridgeman-Bromberg, \cite{BB_advances}) show that such a relation must exist, and work of Epstein-Marden-Markovic \cite{EMM1} and Bridgeman-Canary-Yarmola \cite{BCY} provide sufficient conditions for univalence  in terms of bending norms. 
However, an explicit quantitative bound linking the two descriptions has remained elusive.  

The main contribution of this paper is to provide such an explicit bound. Our approach exploits the geometry of \emph{Epstein surfaces} developed in \cite{Epstein:surface}, a family 
of canonical surfaces in hyperbolic three--space naturally associated to conformal metrics on the disk. By analyzing the convexity and thickness 
properties of these surfaces, we obtain sharp control of convex hull boundaries in terms of the Schwarzian norm. This allows us to deduce explicit 
inequalities relating \(\|Sf\|_\infty\) and \(\|\beta_f\|_L\), valid whenever \(\|Sf\|_\infty\) is sufficiently small.  

As an application, we derive new effective bounds on the bending laminations of quasifuchsian manifolds in terms of the Teichm\"uller distance 
between their conformal boundary components.  

Our results strengthen the analytic--geometric dictionary between the Schwarzian derivative and Teichm\"uller theory. Just as Nehari's estimates 
link small Schwarzian norm to injectivity, our bounds show that small Schwarzian norm also forces controlled geometric bending. In particular, the  use of Epstein surfaces provides a geometric bridge between analytic data on the disk and convex geometric data in hyperbolic three--space.  

\section{Results}
Let  $f:\Delta \rightarrow\hat\CC$ be a locally univalent map of the unit disk. Then  the Schwarzian derivative $Sf$ of $f$ is the holomorphic quadratic differential on the disk given by the formula
$$Sf = \left(\left(\frac{f''}{f'}\right)' -\frac{1}{2} \left(\frac{f''}{f'}\right)^2\right)dz^2$$
The Schwarzian derivative is a measure of how close  $f$ is to being a Mobius transformation. In particular $Sf  = 0$ if and only if $f$ is Mobius  and the solution to $Sf = \phi$ is unique up to post-composition by a Mobius transformation (see \cite{Lehto:book:univalent}).  

 Letting $Q(\Delta)$ be the space of holomorphic quadratic differentials on the unit disk $\Delta$,  then for $\phi \in Q(\Delta)$ we define the pointwise norm by
$$\|\phi(z)\| = \frac{|\phi(z)|}{\rho_h(z)}.$$
where $\rho_h(z) = 4/(1-|z|^2)^2$, and is the hyperbolic area form on $\Delta$. We then define
$$\|\phi\|_\infty = \sup_z \|\phi(z)\|.$$
With this norm we define
$$Q^{\infty}(\Delta)  = \{ \phi \in Q(\Delta)\ | \ \|\phi\|_\infty < \infty\}.$$

We have the following classic result of Nehari.
\begin{theorem}[{Nehari}, \cite{nehari}]
Let $f:\Delta\rightarrow \hat\CC$ be locally univalent. 
If $f$ is univalent then  $\|Sf\|_\infty < 3/2$ and if $\|Sf\|_\infty < 1/2$ then $f$ is univalent.
\end{theorem}

The above plays an important role in Teichm\"uller theory, in particular in Bers description of the complex structure on Teichm\"uller space $\Teich(S)$ and its description as a bounded domain (see \cite{Bers:simunif}).

An alternate description of locally univalent maps is given by Thurston using measured laminations. For a complete description,  see \cite{KT-projective}. His work is more general than we will discuss, giving a parametrization of the space of convex projective structures $\CP(S)$ on a surface $S$  by
$$\CP(S) \simeq \Teich(S)\times \mathcal{ML}(S)$$ 
where $\Teich(S)$ is the space of marked conformal structures on $S$ and $\mathcal{ML}(S)$ is the space of measured laminations on $S$.

We briefly describe Thurston's parametrization in our setting. Given $f:\Delta \rightarrow \hat\CC$ locally univalent, Thurston described a convex hull boundary of the map inside hyperbolic three-space $\Hs$. This  is an immersed locally convex surface in $\Hs$ bent along a  collection of geodesics called bending lines whose bending is described by a transverse measure on the bending lines. The bending lines with this transverse measure gives a measured lamination $\beta_f \in \mathcal{ML}(\Delta)$.

A general measured lamination $\mu \in \mathcal{ML}(\Delta)$ assigns a mass to any arc  $\alpha$ transverse to its support, denoted $i(\mu,\alpha)$. A natural measurement of the size of a measured lamination is the following;
Given an $L > 0$ and $\mu \in \mathcal{ML}(\Delta)$, we define
$$\|\mu\|_L = \sup\{ i(\mu,\alpha)\ | \ \mbox{$\alpha$ open arc transverse to $\mu$ with length $< L$}\}.$$
A measured lamination is {\em uniformly bounded} if $\|\mu\|_L < \infty$ for some  (and hence all) $L>0$ and we define the subset of uniformly bounded measured laminations by $\mathcal{ML}^\infty(\Delta)$.

In this parametrization by measured laminations there are correlate statements to Nehari in terms of $\|\mu\|_L$.  By \cite{BCY} for $L \leq 2\sinh^{-1}(1)$ if $f$ is univalent, then 
$$||\beta_f||_{L} \leq F(L)= 2\cos^{-1}(-\sinh(L/2)).$$
In particular 
$$\mbox{$f$ univalent }\Longrightarrow\ \|\beta_f\|_1 \leq  4.238.$$
Conversely,    Epstein, Marden and Markovic  \cite{EMM1}  proved  that 
 $$||\beta_f||_1 \leq .73 \Longrightarrow \mbox{$f$ univalent.}$$ 
Subsequently  using an approach outlined in unpublished work of Epstein-Jerrard,   \cite{BCY} gave an improved bound by a monotonically increasing function $G:(0,\infty) \rightarrow (0,\pi)$ such that if $||\beta||_{L} < G(L)$ then $f$ is univalent. In particular, this gave
 $$||\beta_f||_1 \leq G(1) = .948 \Longrightarrow \mbox{$f$ univalent.}$$ 
One natural question is, what is the relation between $\|\beta_f\|_L$ and $\|Sf\|_\infty$? 
Combining  Nehari's bounds and the bounds given by $F$ and $G$ we get the explicit relations that
$$ \|Sf\|_\infty \leq \frac{1}{2}  \Longrightarrow \|\beta_f\|_1 \leq   4.238 \qquad\qquad ||\beta_f||_1 \leq  .948\ \Longrightarrow\  \|Sf\|_\infty \leq \frac{3}{2}.$$
Using a compactness argument one also has the following implicit relation.
\begin{theorem}[Bridgeman-Bromberg, \cite{BB_advances}]
Given $L > 0$ there exists a monotonically increasing function $K_L:(0,\infty) \rightarrow (0,\infty)$ such that
if $f:\Delta\rightarrow \hat\CC$ is locally univalent with uniformly bounded bending lamination then
$$\|Sf\|_\infty \leq K_L(\|\beta_f\|_L).$$
\label{BB-bound}
\end{theorem}

In this paper we apply the theory of Epstein surfaces (see \cite{Epstein:surface}) to give the following  explicit bound on bending in terms of the Schwarzian derivative.

\begin{theorem}
Let $f:\Delta\rightarrow \hat\CC$ be univalent with $\|Sf\|_\infty \leq \frac{1}{2}\sech(L)$. Then
$$\|\beta_f\|_L \leq B_L(\|Sf\|_\infty)$$
where
$$B_L(x) = \left\{\begin{matrix}
2\tan^{-1}\left(\frac{2e^{L}x }{\sqrt{1-4x^2}}\right) &  0\leq x \leq \frac{1}{2\sqrt{1+e^{2L}}}\\
\cos^{-1}\left(1-8x^2-4\sinh(L)x\sqrt{1-4x^2}\right) & \frac{1}{2\sqrt{1+e^{2L}}}\leq x \leq \frac{1}{2}\sech(L)
\end{matrix}\right.
$$
Furthermore if $g$ is a univalent map of the complement of $\overline{f(\Delta)}$, then
$$\|\beta_g\|_L \leq B_L(\|Sf\|_\infty).$$
\label{thm:main}\end{theorem}
One application of this is to quasifuchsian manifolds. Given $X, Y \in \Teich(S)$, by Bers simultaneous uniformization (see \cite{Bers:simunif}) there  is an associated quasifuchsian manifold $Q(X,\overline{Y})$ whose conformal boundary is $X\cup\overline{Y}$. This quasifuchsian manifold has an associated convex hull with bending lamination $\beta(X,Y)$.  In this setting we prove the following.

\begin{theorem}
Let $X, Y \in \Teich(S)$ with Teichm\"uller distance $d_T(X,Y) \leq \frac{1}{3}\sech(L)$. Then 
$$\|\beta(X,Y)\|_L \leq B_L\left(\frac{3}{2}d_T(X,Y)\right). $$
\end{theorem}
\medskip
\noindent{\bf The  Function $B_L$:} 
The function $B_L:\left[0,\frac{1}{2}\sech(L)\right]\rightarrow [0,\pi]$ is continuous and monotonically increasing.  Also asymptotically 
$$B_L(x) \simeq 4e^Lx \quad \mbox{for}\quad x \simeq 0.$$
Thus as $\|Sf\|_\infty \rightarrow 0$ then $\|\beta_f\|_L$ is bounded asymptotically linearly by $4e^L\|Sf\|_\infty$. 
{\bf Effectiveness:} To analyse the  effectiveness of this bound, we consider the univalent maps $f(z) = z^c, \ c\in[1,2]$. Then $\|Sf\|_\infty = \frac{1}{2}(c^2-1)$ and  the bending lamination $\beta$ satisfies $\|\beta\|_L = (c-1)\pi$. Thus
$$\|\beta\|_L = \pi\left(\sqrt{2\|Sf\|_\infty+1}-1\right) \simeq \pi\|Sf\|_\infty \quad \mbox{for}\quad\|Sf\|_\infty \simeq 0.$$
Thus the asymptotic comparison is between the functions $\pi x$ and $4e^Lx$ which is close for $L$ bounded.
 
{\bf  $L = 1$:} A natural choice of parameter is $L = 1$. Then $B_1$ gives the function $B_1:  [0, .324] \rightarrow [0,\pi]$ graphed below.

\begin{figure}[htbp] 
   \centering
   \includegraphics[width=2in]{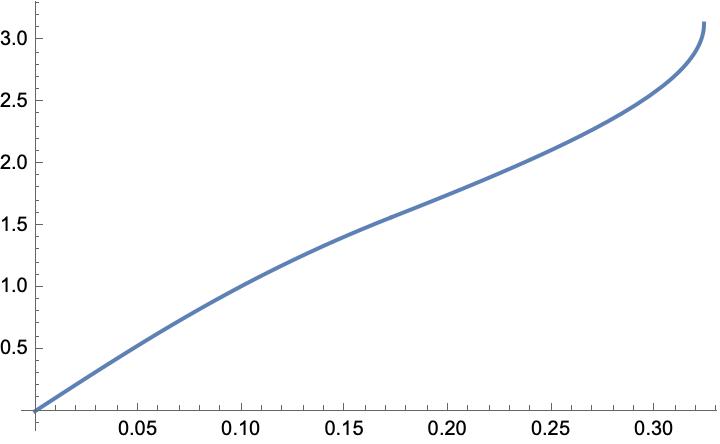} 
   \caption{Graph of $B_1$}
   \label{fig:bend_graph}
\end{figure}

\section{Background}

\subsection{Teichmuller space, hyperbolic 3-manifolds, convex hulls}
Given $S$ a closed surface, the {\em Teichm\"uller space} of $S$, denoted $\Teich(S)$, is the space of marked conformal structures on $S$. Specifically 
 $$\Teich(S) =\{ (f:S\rightarrow X)\ | \mbox{ $f$ is  a diffeomorphism, $X$ a Riemann surface}\}/\sim$$
where $(f:S\rightarrow X) \sim  (g:S\rightarrow Y)$ if $g\circ f^{-1}:X\rightarrow Y$ is homotopic to a conformal map.

By Riemann uniformization, $\Teich(S)$ is also the space of marked hyperbolic structures on $S$.

We now consider hyperbolic 3-manifolds. A complete hyperbolic 3-manifold is a quotient manifold $M = \Hs/\Gamma$ where  $\Gamma$ is a discrete subgroup of $\psl \simeq Isom^+(\Hs)$. The {\em limit set} $\Lambda(\Gamma)$ of $\Gamma$ is  defined to be
$$\Lambda(\Gamma) = \overline{\Gamma x}\cap \hat\CC$$
where $x \in \Hs$ is any point. The {\em domain of discontinuity} is $\Omega_\Gamma = \hat\CC-\Lambda(\Gamma)$ and the {\em conformal structure at infinity} is
$$\partial_c M = \Omega_\Gamma/\Gamma.$$
The {\em convex hull} $H(\Lambda(\Gamma))$ is the smallest convex subset of $\Hs$ containing all the geodesics with both endpoints in $\Lambda(\Gamma)$.
The {\em convex core} is 
$$C(M) = H(\Lambda(\Gamma))/\Gamma.$$
 Thurston showed that the components of the convex hull boundary are given by convex pleated planes (see \cite{Thurston:book:GTTM}). That is, for each component $C$ of $\partial H(\Lambda(\Gamma))$, there is  a measured lamination $\mu$ on $\Hp$ and a homeomorphism $f:\Hp \rightarrow C \subseteq \Hs$ such that $f$ is  an isometry in the complement of the support of $\mu$ and the bending of $f$ is along the support of $\mu$ given by the transverse measure on $\mu$. 
  
 To be more precise, we describe the transverse measure $\mu$. For complete details, see  Epstein-Marden's paper \cite{Epstein:Marden:convex}.  If $x \in \partial H(\Lambda(\Gamma))$, a {\em support half-space} $H$ to $x$ is a half-space $H$ whose interior is disjoint from $H(\Lambda(\Gamma))$ and $x \in \partial H$.  We let $m$ be the collection of bending lines of $C$. For $\alpha:[0,1]\rightarrow C$ an arc transverse to $m$, given a partition $\mathcal P = \{0 = t_0 < t_1 < \ldots t_n =1\}$ and support half-spaces $\mathcal H =\{ H_{t_i} \ | \ \alpha(t_i)\in \partial H_{t_i}\}$ then we let
 $$i_{(\mathcal P,\mathcal H)}(\alpha,\mu) = \sum_{i=0}^{n-1} ext\angle H_{t_i},H_{t_{i+1}}$$
 where $ext\angle H_a, H_b$ is the exterior angle between half-spaces $H_a$ and $H_b$.
 Then the transverse measure on $\alpha$ is given by
 $$i(\alpha,\mu) = \lim_{|\mathcal P|\rightarrow 0} i_{(\mathcal P,\mathcal H)}(\alpha,\mu).$$
In our work, we will need to bound the transverse measure. We call the pair $(\mathcal P, \mathcal H)$  {\em good} if for any $H_t$ support half-space to $\alpha(t)$ with $t \in [t_i, t_{i+1}]$ then $H_t$ intersects both  $H_{t_i}, H_{t_{i+1}}$. If  $(\mathcal P, \mathcal H)$ is good then by elementary hyperbolic geometry $i_{(\mathcal P,\mathcal H)}(\alpha,\mu)$ is monotonically decreasing under refinement. Then the transverse measure on $\alpha$ satisfies
 $$i(\alpha,\mu) = \inf\{ i_{(\mathcal P,\mathcal H)}(\alpha,\mu)\ | \ (\mathcal P, \mathcal H) \mbox{ good}\}.$$
 In particular for  $(\mathcal P, \mathcal H)$ good, one useful bound is
 \begin{equation}
 i(\alpha,\mu) \leq i_{(\mathcal P,\mathcal H)}(\alpha,\mu).
 \label{good}
 \end{equation}
 
One type of  hyperbolic 3-manifold closely related to the theory of univalent maps are quasifuchsian manifolds. A {\em quasifuchsian group} $\Gamma$ is a Kleinian group with limit set  a Jordan curve and whose action preserves each component of its complement. Then by Thurston, the convex hull boundary is the union of two convex pleated planes (see \cite{Thurston:book:GTTM}). By Bers simultaneous uniformization (see \cite{Bers:simunif}), if $X, Y$ are conformal structures on a closed surface $S$, then there exists quasifuchsian group $\Gamma$ with conformal boundary at infinity $X\cup \overline{Y}$. Bers showed further that this simultaneous uniformization gives a homeomorphism between the space of quasifuchsian structures $QF(S)$ on a closed surface $S$  and $\Teich(S)\times\Teich(\overline{S})$.

 \subsection{Thurston's parametrization for locally valent maps}

We define 
$$P(\Delta) =\{f:\Delta\rightarrow \hat\CC\ | \ \mbox{ $f$ locally univalent}\}/\sim$$
where $f\sim g$ if $g = m\circ f$ for $m \in \psl$. This can be identified as the space of complex projective structures on $\Delta$ (see \cite{KT-projective}). Taking the Schwarzian we can identify $P(\Delta) = Q(\Delta)$ the space of holomorphic quadratic differentials. 

On $Q(\Delta)$ we define the pointwise norm for $\phi \in Q(\Delta)$ by
$$\|\phi(z)\|= \frac{|\phi(z)|}{\rho_h(z)}$$
where $\rho_h(z) = 4/(1-|z|^2)^2$ is the hyperbolic metric on $\Delta$.
We define the subspace
$$Q^\infty(\Delta) =\{  \phi\ |\  \|\phi\|_\infty < \infty\}.$$ 

We now describe Thurston's parametrization of $P(\Delta)$ by measured laminations. See \cite{KT-projective} for further details.

We take the approach of Bonahon in describing the space of measured laminations (see \cite{bonahon:currents}). We let $G(\Hp)$ be the space of unoriented geodesics. Then identifying the boundary of $\Hp$ with $\mathbb S^1$ then $G(\Hp) \simeq (\mathbb S^1 \times \mathbb S^1)-\mbox{diag})/\Z_2$ where $Z_2$ acts by $(x,y) \rightarrow (y,x)$. A {\em geodesic lamination} is a closed subset of $G(\Hp)$ whose points are mutually disjoint as geodesics. A {\em measured lamination} on $\Hp$ is a Borel measure on $G(\Hp)$ whose support is a geodesic lamination. The space of measured laminations on $\Hp$ is denoted $\mathcal{ML}(\Hp)$ and given the weak$^*$ topology.

Given $\mu \in \mathcal{ML}(\Hp)$, and $\alpha$ an arc transverse to the support of $\mu$, the {\em transverse measure} on $\alpha$ is defined to be
$$i(\mu,\alpha) = \mu(G(\alpha))$$
where $G(\alpha)$ is the set of geodesics intersecting $\alpha$ transversely.  For $L>0$ we define

$$\|\mu\|_L = \sup\{i(\mu,\alpha)\ | \ \mbox{ $\alpha$ open transverse to $\mu$ of length $< L$}\}.$$
Then we define the set of {\em uniformly bounded} measured laminations
$$\mathcal{ML}^\infty(\Delta)\ =   \{ \mu\ |\  \|\mu\|_L < \infty \mbox{ for some $L$}\}.$$

We first describe Thurston's parametrization  for  univalent maps. If $f:\Delta \rightarrow\hat\CC$ is univalent with $f(\Delta) = \Omega_f$ then given an open round disk $D \subseteq \Omega_f$ we let $H_D$ be the  half-space in $\Hs$ with boundary $D$. If $D$ is maximal, then $H_D$ is called a support half-space.
Then we define the dome of $\Omega_f$ by
$$ \dome(f) = \bigcap_{D \mbox{ maximal}} (H^o_D)^c$$
By definition $\dome(f)$ is closed and convex. The $\dome(f)$ is also equal to the  convex hull of the complement of $\Omega$.  

By work of Thurston $\partial \dome(f)$ is topologically a disk and has intrinsic metric, the hyperbolic metric. Thus $\partial Dome(f)$ is isometric to $\Hp$. Furthermore there is an isometry map $F:\Hp \rightarrow \partial \dome(f)$ which is isometric in the complement of the support of measured lamination $\mu$ with measure given by the bending of $\partial \dome(f)$ along a geodesic lamination $m$. Thurston's parametrization of $[f] \in P(\Delta)$ is this measured lamination $\mu \in \mathcal{ML}(\Hp)$.
 
 Although we will only be considering Thurston's parametrization for univalent maps, we briefly describe the parametrization for the general (locally-univalent) case. We let  $f :\Delta\rightarrow \hat\CC$ be a locally univalent map. We first define a round disk for $f$ to be an open disk in $U \subseteq \Delta$  such that $f:U \rightarrow f(U)$ is a univalent map where $f(U)$ is a round disk in $\hat\CC$. Given $f:\Delta\rightarrow \hat\CC$ we consider 
$$\mathcal U_f = \{ U\ | \ \mbox{ $U$ is a maximal round disk for $f$}\}$$
For each maximal disk $U$ the image $f(U)$ is a round disk and is the boundary of a unique halfspace $H_{f(U)}$ in $\Hs$. We then define the $\dome(f)$ as before and now obtain a map $F:\Hp \rightarrow \Hs$ and a measured lamination $\mu$ which is an isometry on each component of the complement of $\mu$ and has bending given as above. Unlike the univalent case, the map $F$ is not a homeomorphism but we still obtain a measured lamination.
 
 By Thurston the above gives a homeomorphism  $\Psi:P(\Delta) \rightarrow \mathcal{ML}(\Hp)/\sim$ where $\mu\sim \nu$ if $\nu = m^*\mu$ for $m \in \pslr$.

\section{Bounding bending by thickness of convex hull}
For  a Jordan curve  $\gamma \subseteq \hat\CC$ we  define the convex hull $H(\gamma)$  in $\Hs$ as before as the smallest convex set in $\Hs$ containing all the geodesics with both endpoints in $\gamma$.  
Such hulls arise as the convex hulls of quasifuchsian groups. If $\gamma$ has complement given by Jordan domains $\Omega_1,\Omega_2$ then it is easy to see that 
$$H(\gamma) = \dome(\Omega_1)\cap \dome(\Omega_2).$$
We label the boundary components of $H(\gamma)$ by   $C_1, C_2$ where $C_i = \partial \dome(\Omega_i)$. Furthermore $C_i$ is a convex pleated plane with bending lamination $\beta_i$.
 We define the {\em thickness} of $H(\gamma)$ 
$$T_1(\gamma) = \sup \{d(x, C_2)\ | \ x\in C_1\} \qquad T_2(\gamma) = \sup \{d(x, C_1)\ | \ x\in C_2\}.$$
This could be infinite but for $\gamma$ equal the limit set of a convex cocompact quasifuchsian group, it is always finite. 

We define the following function $C_L$.
$$C_L(r) =\left\{ \begin{matrix}
2\tan^{-1}(e^L\sinh(r))& 0 \leq \sinh(r) \leq e^{-L}\\
\cos^{-1}\left(1-2\tanh^2(r)\left(1+\frac{\sinh(L)}{\sinh(r)}\right)\right) &  e^{-L} \leq \sinh(r) \leq 1/\sinh(L) \end{matrix}
\right.$$
We prove the following bound using elementary hyperbolic geometry.

\begin{theorem}
Let $L > 0 $  and $\gamma$ a Jordan curve. Let $\alpha: [0,L]\rightarrow C_1$ be a geodesic parameterized by arc length with $\alpha(0) = x$. If $\sinh(L)\sinh(d(x,C_2)) \leq 1$ then
$$i(\alpha,\beta_1) \leq  C_L(d(x,C_2)) $$
In particular if $\sinh(L)\sinh(T_i(\gamma)) \leq 1$ then
$$\|\beta_i\|_L \leq  C_L(T_i(\gamma)) $$
\end{theorem}

{\bf Proof:} 
 We let $H_1$ be a support halfspace to $x_1 = \alpha(0)$ and $H_2$ a support halfspace to $x_2 = \alpha(L)$. Let $x_3 \in C_2$ with $d(x_1, x_3) = d(x, C_2)$. We choose a support halfspace $H_3$ at $x_3$ to $C_2$. By definition $H_3$ is disjoint from $H_1$ and $H_2$. Let $\partial H_i$ be the boundary planes. We choose $H$ to be the unique plane perpendicular to all $\partial H_i$. Projecting perpendicularly $H_i$ project to halfplanes $H'_i$, with $H'_3$ disjoint from $ H'_1$ and $H'_2$. The points $x_i$ project to points $x'_i$ and as perpendicular projection  is distance non-increasing, $d(x'_1,x'_2) \leq L$ and $d(x'_1,x'_3) \leq d(x_1, C_2)$. Further $x'_1 \not\in (H'_2)^o$ and $x'_2 \not\in (H'_1)^o$. Also the exterior angle between $H_1$ and $H_2$ is the exterior angle between $H'_1, H'_2$.   Thus by Lemma \ref{lemma:tech}  below if $\sinh(L)\sinh(d(x_1,C_2)) \leq 1$ then $H'_1, H'_2$ intersect with 
$$ext\angle H_1, H_2 = ext\angle H'_1, H'_2 \leq C_L(d(x_1,C_2)).$$
Therefore as the pair of support planes $H_1, H_2$ give a  good pair $(\mathcal P,\mathcal H)$ with $
\mathcal P = \{0,1\}, \mathcal H = \{ H_1, H_2\}$ (see equation \ref{good}) then
$$i(\alpha,\beta_1) \leq i_{(\mathcal P, \mathcal H)}(\alpha,\beta_1)   \leq C_L(d(x_1,C_2)).$$
The bound on $\|\beta_i\|_L$ in terms of thickness follows by definition of the norm.
 \eproof

\subsection{Hyperbolic Trigonometry}
We will make use of the following hyperbolic trigonometry formulae for a triangle with one ideal vertex . Let $T$ be a hyperbolic triangle with angles $\alpha,\beta,\gamma$ and sides $A,B,C$. If $\gamma = 0$ then
$$\cosh(C) = \frac{\cos(\alpha)\cos(\beta)+1}{\sin(\alpha)\sin(\beta)}$$ and
$$\sinh(C) = \frac{\cos(\alpha)+\cos(\beta)}{\sin(\alpha)\sin(\beta)}\qquad \tan(\alpha/2)\tan(\beta/2) = e^{-C}.$$
The first is the standard hyperbolic cosine formula (see \cite{Thurston:book:GTTM}) and the other two we could not find a reference for but can be easily derived from the first.
We  will need the following elementary lemma involving half-planes.

\begin{lemma}
Let $H_1, H_2, H_3$ be half-planes in $\Hp$ with $H_3$ disjoint from $H_1$ and $H_2$. Further let $z_i \in \partial H_i$ be points such that $d(z_1, z_2) \leq L$ and  $d(z_1,z_3) \leq r$ with $z_1 \not\in H^o_2, z_2\not\in H_1^o$.
If $\sinh(r)\sinh(L) \leq 1$ then  $H_1, H_2$ intersect with exterior angle $\theta \leq C_L(r).$
\label{lemma:tech}
\end{lemma}

{\bf Proof:}
We place $z_1$ at the origin in the Poincare model and let $H_1 =\{ z \in \Delta \ | \Im(z) < 0\}$.
Let $g_i$ be the geodesic boundary of $H_i$ and $\phi_i$ be visual angle of $H_i$ from $x_1$. Then $\phi_1 = \pi$ and for $\phi_i, \ i \neq 1$, we have a triangle with angles $\phi_i/2, \pi/2, 0$ and side length $d(x_1, H_i)$.  We define angles $\hat\phi_2,\hat\phi_3 \in (0,\pi)$ as follows; 
$$\tan(\hat\phi_2/2) = \frac{1}{\sinh(L)}  \qquad \tan(\hat\phi_3/2) = \frac{1}{\sinh(r)}$$ 
As $d(z_1, H_2) \leq d(z_1,x_2) \leq L$ then by the above trigonometry formulae
$$\tan(\phi_2/2) = \frac{1}{\sinh(d(x_1, H_2))} \geq \frac{1}{\sinh(L)} = \tan(\hat\phi_2/2)$$
Similarly 
$$\tan(\phi_3/2) \geq \frac{1}{\sinh(r)} = \tan(\hat\phi_3/2).$$
In particular $ \hat\phi_i  \leq\phi_i $ for $i = 2,3$. As $\sinh(r)\sinh(L) \leq 1$ then 
$$\tan(\phi_2/2)\tan(\phi_3/2) \geq 1.$$
As
$$\tan((\phi_2+\phi_3)/2) = \frac{\tan(\phi_2/2)+\tan(\phi_3/2)}{1-\tan(\phi_2/2)\tan(\phi_3/2)}$$
it follows that $(\phi_2+\phi_3)/2 \geq \pi/2$ and $\phi_2+\phi_3 \geq \pi$.
Thus the total angle subtended by $H_1,H_2$ and $H_3$ is greater than $2\pi$. It follows that $H_1, H_2$ intersect.

We now move to proving the bound.  If $H_1 \subseteq H_2$ then as $z_1$ is not in the interior of $H_2$  then $z_1 \in g_1\cap g_2$. It follows that $H_1 = H_2$  and the exterior angle is zero and the result holds. Similarly for $H_2 \subseteq H_1$. 
 
From the above, if $H_1 \neq H_2$ then $g_1,g_2$ intersect transversely with $g_1\cap g_2 = t \in (-1,1) \subset \R$. As $g_1, g_2$ intersect transversely, we can assume that $1$ is in the boundary of $H_2$. 
Thus it follows that $t \geq 0$ as otherwise $z_1 \in H^o_2$.

We let $g_2$ have endpoints $p = e^{ia}, q = e^{ib}$ where $0 < a < \pi$ and $\pi < b  < 2\pi$. As $z_2 \not\in H_1^o$ then $z_2$ is on the geodesic ray  $\overrightarrow{tp}$. Thus $\overrightarrow{tp}$ intersects the ball of radius $L$ about $z_1$. 

We let $h$ be the geodesic perpendicular to $g_1$  at the point a distance $L$ from $z_1$ on the positive real axis. Then by definition of $\hat\phi_2$, $h$ has endpoints $e^{i\hat\phi_2/2}, e^{-i\hat\phi_2/2}$.

\begin{figure}[htbp] 
   \centering
   \includegraphics[width=3in]{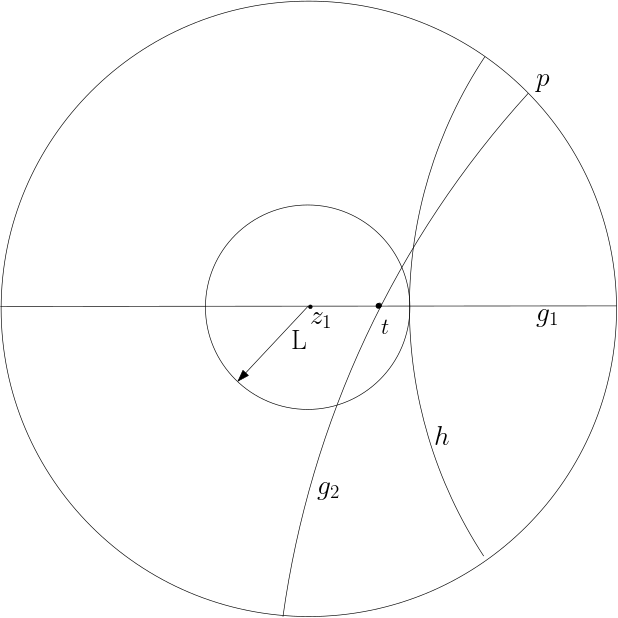} 
   \caption{Case 1 configuration of geodesics}
   \label{figure:triangle}
\end{figure}

It follow that if $a \leq \hat\phi_2/2$ then $d(z_1,t) \leq L$ as otherwise the ray $\overrightarrow{tp}$ does not intersect the ball of radius $L$ about $z_1$ (see figure \ref{figure:triangle}).

We consider two cases.

{\bf Case 1, $a \leq \hat\phi_2/2$:}   As $a\leq \hat\phi_2/2$ then $d(z_1, t) \leq L$ and $g_2$ meets the positive real axis in angle $\theta \leq \pi/2$. We take the triangle $z_1, t, p$ labelling the angles at $z_1, t$ by $\alpha,\beta$ respectively and the side $C = d(z_1, t)$. Then 
$$e^{-L} \leq e^{-C} = \tan(\alpha/2)\tan(\beta/2) .$$
We note that $\alpha = a$ and $\beta = \pi-\theta$. Thus
$$e^{-L} \leq \tan(\alpha/2)\tan(\beta/2) = \tan\left(\frac{a}{2}\right)\cotan\left(\frac{\theta}{2}\right) .$$
Thus
$$\tan(\theta/2) \leq e^{L}\tan(a/2)$$

{\bf Case 2, $a \geq \hat\phi_2/2$:}
 As $g_3$ must intersect the ball of radius $L$ about $z_1$, we consider the geodesic $k$ with endpoint $p=e^{ia}$ and tangent to the circle of radius $L$ about $z_1$. For $a \leq \hat\phi_2$ then $k$ intersects the positive $x$-axis. Let $\phi$ be the angle $k$ intersects the positive $x$-axis. It follows that $\theta \leq \phi$. We now observe taking the perpendicular from $z_1$ to $k$ there is a right-angled hyperbolic triangle with side of length $L$ opposite angle $\pi-\phi$ and  angle equal $a- \hat\phi_2/2$ at $z_1$. Then
 $$\cosh(L) = \frac{\cos(\pi-\phi)}{\sin(a-\hat\phi_2/2)}.$$
 Thus
 \begin{align*}
 \cos(\phi) &= -\cos(\pi-\phi) = -\cosh(L) \sin(a-\hat\phi_2/2) \\
 &= -\cosh(L)\left(\sin(a)\cos(\hat\phi_2/2)-\cos(a)\sin(\hat\phi_2/2)\right)\\
 &= -\cosh(L)\left(\sin(a)\tanh(L)-\cos(a)\sech(L)\right)\\
 &= -\sin(a)\sinh(L)+\cos(a)
 \end{align*}
 
 giving
 $$\theta \leq \phi = \cos^{-1}(-\sin(a)\sinh(L)+\cos(a)).$$
 
  \begin{figure}[htbp] 
   \centering
   \includegraphics[width=3in]{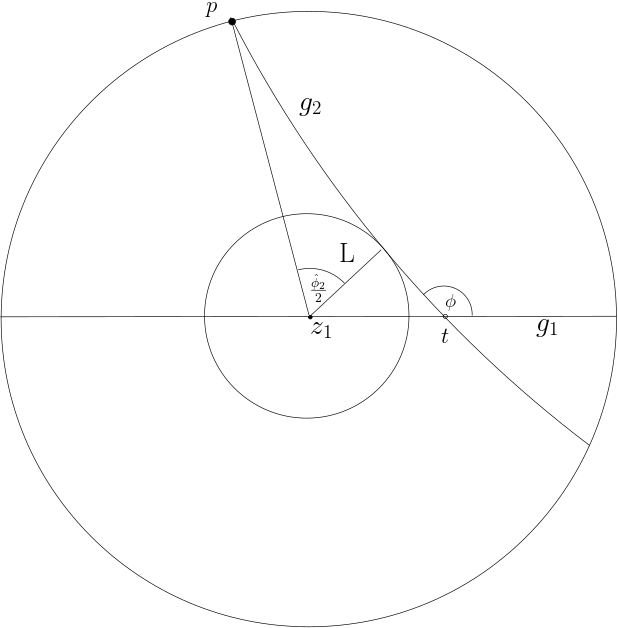} 
   \caption{Case 2 configuration of geodesics}
   \label{figure:triangle2}
\end{figure}

  Thus we have the bound $\theta \leq f(a)$ where function $f$ is
$$ f(a) = \left\{\begin{matrix}
2\tan^{-1}(e^L\tan(a/2)) & 0 \leq a \leq \hat\phi_2/2\\
\cos^{-1}(-\sin(a)\sinh(L)+\cos(a)) &  \hat\phi_2/2  \leq a \leq \hat\phi_2
\end{matrix}\right.
 $$
 We note that $f:[0, \hat\phi_2] \rightarrow [0,\pi]$, is continuous and strictly monotonically increasing. Also  interval $[0,\hat\phi/2]$ is mapped to $[0,\pi/2]$ and interval $[\hat\phi/2,\hat\phi]$ is mapped to $[\pi/2,\pi]$.
 
Now we obtain a bound in terms of $r, L$. As $H_2,H_3$ are disjoint then 
$$a \leq \pi-\phi_3 \leq \pi-\hat\phi_3.$$ 
We note that as $\sinh(r)\sinh(L) \leq 1$ then $\hat\phi_2+\hat\phi_3 \geq \pi$ giving $\pi-\hat\phi_3 \leq \hat\phi_2$ and $\pi-\hat\phi_3$ is in the domain of $f$. Thus we can define for $\sinh(r)\sinh(L) \leq 1$ the function $C_L(r)$ by
$$\theta \leq f(a) \leq f(\pi-\hat\phi_3) = C_L(r)$$

To check this gives the desired formula, we first note that
$$\tan((\pi-\hat\phi_3)/2) = \cotan(\hat\phi_3/2) =\sinh(r).$$
Thus for $\pi-\hat\phi_3 \leq \hat\phi_2/2$ then
$$C_L(r) =  f(\pi-\hat\phi_3) = 2\tan^{-1}(e^L\sinh(r)).$$
We now confirm the formula  for $\pi-\hat\phi_3 \geq \hat\phi_2/2$.

 \begin{align*}
 C_L(r) &= f(\pi-\hat\phi_3) \\
 &=  \cos^{-1}\left(-\sin(\pi-\hat\phi_3)\sinh(L)+\cos(\pi-\hat\phi_3)\right)\\
 &=  \cos^{-1}\left(-\sin(\hat\phi_3)\sinh(L)-\cos(\hat\phi_3)\right)\\
 &=  \cos^{-1}\left(-2\sin(\hat\phi_3/2)\cos(\hat\phi_3/2)\sinh(L)-(2\cos(\hat\phi_3/2)-1)\right)\\
  & = \cos^{-1}\left(1-2\tanh^2(r) -\frac{2\tanh(r)\sinh(L)}{\cosh(r)}\right)\\
  & = \cos^{-1}\left(1-2\tanh^2(r)\left(1+\frac{\sinh(L)}{\sinh(r)}\right)\right)
 \end{align*}
 as desired. 
 
 Finally we note that in terms of $L,r$ the equation $\pi-\hat\phi_3 = \hat\phi_2/2$ gives
$$e^{-L} = \tan(\hat\phi_2/4) = \tan((\pi-\hat\phi_3)/2) = \sinh(L).$$
Thus the piecewise intervals are $0\leq \sinh(r) \leq e^{-L}$ and $e^{-L}\leq \sinh(r) \leq 1/\sinh(L)$.
 \eproof

\section{Schwarzian bound on Thickness}

We will first use Epstein surfaces to bound the thickness of the convex hull. 
We have some notation. If $f:\Delta \rightarrow \Hs$ is an immersion onto surface $S$, we define the {\em fundamental pair} $(g, B)$ to be the pullback of the induced metric and the shape operator on $f(\Delta)$ pulled back to $\Delta$ respectively. The eigenvalues of $B$ are the {\em principal curvatures} of the surface $S$.

If $B$ does not have eigenvalues equal to $-1$ we let $$\hat g = (\id+B)^* g\qquad \hat B = (\id-B)(\id+B)^{-1}.$$ Then $(\hat g,\hat B)$ is called the {\em dual pair} for $(g,B)$. One reason to consider this dual pair is that they are an equivalent representation of the pair $(g,B)$ and often have a simpler description. We observe that $B$ has eigenvalues in $[0,\infty)$ if and only if $\hat B$ has eigenvalues in $(-1,1]$. 

In \cite{Epstein:surface}, C. Epstein showed how to associate to a conformal metric on the domain in $\hat\CC$ an immersed surface in $\Hs$, called the {\em Epstein surface} of the conformal metric. If one takes the hyperbolic metric on the domain of discontinuity, then this surface is called the {\em Poincar\'e-Epstein} surface. 

We have the following properties of the Poincar\'e Epstein surface.

\begin{theorem}[{Epstein, \cite{Epstein:surface}}]
Let $f:\Delta \rightarrow \hat\CC$ be a  univalent map with $\|Sf\|_\infty < 1/2$. Then the Poincare-Epstein surface for $f$ is an embedded surface $Ep_f:\Delta \rightarrow \Hs$ with principal curvatures at $Ep_f(z)$ equal
$$-\frac{\|\phi(z)\|}{\|\phi(z)\|\pm1}.$$
Furthermore normal flow on the surface gives a foliation of $\Hs$.  If  $Ep^t_f:\Delta \rightarrow \Hs$  is the surface given by time $t$ normal flow, then $\lim_{t\rightarrow \infty}Ep^t_f(z) = f(z)$.
\end{theorem}
We now list some properties of the Poincare-Epstein surface of $f$ that follow directly from the definition of $\hat B$ (see \cite{BB-OS} for more details).

\begin{itemize}
 \item The dual shape operator $\hat B$ has eigenvalues $1\pm 2\|\phi(z)\|$ at $z$.
 \item The surface $Ep^t_f:\Delta \rightarrow \Hs$   has dual shape operator $\hat B_t = e^{-2t}\hat B$. 
 \item If $e^{2t} \geq 1+2\|\phi\|_\infty$ the surface $Ep^t_f$ is locally convex (principal curvatures both non-negative).
 \item If  $e^{2t} < 1-2\|\phi\|_\infty$  the surface $Ep^t_f$ is locally concave (principal curvatures both non-positive).
 \end{itemize}

Using the Poincar\'e-Epstein surface we obtain the following bound on thickness.
\begin{corollary}
Let $f:\Delta \rightarrow \Omega$ be univalent with $\|Sf\|_\infty < 1/2$ with   $\partial\Omega = \gamma$. Then
$$T_i(\gamma) \leq \frac{1}{2}\log\left(\frac{1+2\|Sf\|_\infty}{1-2\|Sf\|_\infty}\right).$$
\end{corollary}

{\bf Proof:}
We let $e^{2t_0} > 1+2\|Sf\|_\infty$ and $e^{2t_1} < 1-2\|Sf\|_\infty$. Let 
$$N(t_0,t_1) = \{Ep_f^t(z)\ | \ z\in \Delta, \ t \in[t_0,t_1]\}.$$ As normal flow gives a foliation for $\Hs$, $N$ is foliated by disjoint geodesic arcs of length $t_0-t_1$. Also by convexity/concavity of surfaces $Ep_f^{t_0}, Ep_f^{t_1}$ then $N$ is convex and by minimality $H(\gamma) \subseteq N$. Thus $z \in  H(\gamma)$ is on a line segment of length $t_0-t_1$ connecting the Epstein surfaces $Ep_f^{t_0}$ to $Ep_f^{t_1}$.  In particular, every point of $\partial H(\gamma)$ is on an arc of length $t_0-t_1$ containing a point of the other boundary component of $\partial H(\gamma)$. Thus
$$T_i(\gamma) \leq t_0-t_1.$$
As we can choose $t_i$ such that $e^{2t_i}$ are arbitrarily close to $1+2\|Sf\|_\infty, 1-2\|Sf\|_\infty$, the result follows.
\eproof

We now prove Theorem \ref{thm:main}. 

{\bf Proof of Theorem  \ref{thm:main}:}
We let
$$r(s) = \frac{1}{2}\log\left(\frac{1+2s}{1-2s}\right).$$
and define
$$B_L(s) = C_L(r(s)).$$
We have 
$$\tanh(r(s)) = 2s \qquad \sinh(r(s)) = \frac{2s}{\sqrt{1-4s^2}}.$$
Thus if $\sinh(r(s)) = a$ then
$$s = \frac{1}{2}\tanh(r(s)) = \frac{1}{2}\frac{\sinh(r(s))}{\cosh(r(s))} = \frac{1}{2}\frac{a}{\sqrt{1+a^2}}.$$
Thus the domain $\sinh(L)\sinh(r(s)) \leq 1$ corresponds to
$$s \leq \frac{1/\sinh(L)}{2\sqrt{1+1/\sinh^2(L)}} = \frac{1}{2}\sech(L)$$
and the domain $e^{L}\sinh(r(s)) \leq 1$ corresponds to
$$s \leq \frac{1}{2\sqrt{1+e^{2L}}}.$$
By the above
$$T_i(f) \leq r(\|Sf\|_\infty).$$ 
Therefore by monotonicity of $C_L$, then for $\|Sf\|_\infty < \frac{1}{2}\sech(L)$, we have
$$\|\beta\|_L \leq C_L(T_i(f)) \leq C_L(r(\|Sf\|_\infty)) = B_L(\|Sf\|_\infty).$$
where
$$B_L(x) = \left\{\begin{matrix}
2\tan^{-1}\left(\frac{2e^{L}x }{\sqrt{1-4x^2}}\right) &  0\leq x \leq \frac{1}{2\sqrt{1+e^{2L}}}\\
\cos^{-1}\left(1-8x^2-4\sinh(L)x\sqrt{1-4x^2}\right) & \frac{1}{2\sqrt{1+e^{2L}}}\leq x \leq \frac{1}{2}\sech(L)
\end{matrix}\right.
$$
\eproof

\section{Teichm\"uller distance}
The bound on bending in terms of the Teichm\"uller distance will follow by bounding the derivative of the Bers map.

Given $X \in \Teich(S)$, we define a map $\Phi_X:\Teich(\overline{S}) \rightarrow Q(X)$ where $\Phi_X(Y)$ is the Schwarzian derivative of the map uniformizing the domain corresponding to $X$ in the quasifuchsian manifold with conformal boundary $X\cup Y$. By Ahlfors-Weill we have the following.

\begin{theorem}[{Ahlfors-Weill \cite{ahlforsweill}}]
Let $\|\Phi_X(Y)\|_\infty < 1/2$ then
$$d_T(X,\overline{Y})\leq \tanh^{-1}(2\|\Phi_X(Y)\|_\infty).$$
\end{theorem}

In \cite{TT-WP} Takhtajan and Teo consider the Lipschitz constant for the Bers mapping  and show that it is $12$-Lipschitz with respect to the $L^2$-metric on both domain and range. Modifying their proof by using the Area theorem, we can improve this to $3/2$-Lipschitz for both the $L^2$ and $L^\infty$ metrics.

\begin{theorem}
The map $\Phi_X$ is $3/2$-Lipschitz with respect to the Teichm\"uller metric on $\Teich(S)$ and the $L^\infty$ norm on $Q(X)$. In particular
$$\|\Phi_X(Y)\|_\infty \leq \frac{3}{2}d_T(X,\overline{Y}).$$
\end{theorem}

The bound on bending in terms of the Teichm\"uller distance follows immediately. 
In order to prove the Lipschitz bound, we will need to consider the integral formula for the derivative of the Bers embedding using the complex analytic structure on Teichm\"uller space. For full details see  Imayoshi and Taniguchi's book \cite{Imayoshi:Taniguchi:book}.

We let $\Half = \{ z\ | \ \Im(z) > 0\}$ be the upper half-plane. 
For $X = \Half/\Gamma$ we define $B(\Half,\Gamma)$ to be the set of $\Gamma$ invariant beltrami differentials and $Q(\Half,\Gamma)$ be the space of holomorphic quadratic differentials on $\Hp$ invariant under $\Gamma$. Then for $\mu \in B(\Half,\Gamma)_1$, the open unit ball in $B(\Half,\Gamma)$, we let  $\hat\mu$ be the Beltrami differential on $\hat\CC$ given by
$$\hat\mu(z) = 
\left\{\begin{matrix} 
\mu(z), & z\in \Half\\
&\\
\overline{\mu(\overline{z})} & z\in \overline{\Half}
\end{matrix}\right.
$$
Then we define $f_\mu:\Half\rightarrow \Half$ to be the restriction to $\Half$ of the unique solution to the Beltrami equation $F_{\overline z} = \hat\mu F_{z}$, fixing $0,1,\infty$. Then we can identify 
$$\Teich(S) = B(\Half,\Gamma)_1/\sim$$
where the equivalence relation $\mu \sim \nu$ if $f_\mu$ and $f_\nu$ are equal on $\overline{\R}$. 
Then we have
$$T_X\Teich(S) \simeq B(\Half,\Gamma)/N(\Gamma)$$
where 
$$N(\Gamma) = \left\{ \mu\ \left| \ \int_X \mu \phi = 0 \quad \forall \quad \phi \in Q(\Half,\Gamma)\right\}\right..$$
The $L^p$ norm on $T_X\Teich(S)$ is given by
$$\|[\mu] \|^p_p = \inf_{\mu \in [\mu]} \int_{\Half/\Gamma} |\mu(z)|^p \rho(z)|dz|^2 .$$
Then for $\phi \in Q(\Half,\Gamma)$ we define the pointwise norm by $\|\phi(z)\| = |\phi_h(z)|/\rho(z)$ and the $L^p$ norm by
$$\|\phi\|^p_p = \int_{\Half/\Gamma} \|\phi(z)\|^p\rho_h(z) |dz|^2.$$

Given $Y \in \Teich(S)$ with $Y = [\mu]$ then we have quasi-conformal map $f_{\mu}:\CC \rightarrow \CC$ which has Beltrami differential $\mu$ on $\Half$ and $0$ on $\overline{\Half}$. The Bers embedding then lifts to the map $\Phi: B(\Half,\Gamma)_1 \rightarrow Q(\overline{\Half},\Gamma)$ given by
$$\Phi(\mu) = S(f_\mu).$$

If $\mu \in  B(\Half,\Gamma)_1$ and $\nu \in B(\Half,\Gamma)$ then letting
$$\mu_t = \mu+t\nu +O(t^2) \in B(\Half,\Gamma)_1$$ we obtain a deformation of $Y$ given by $f_{\mu_t}$. We define the derivative by
$$\dot\Phi_{\mu}([\nu])(z) = \lim_{t\rightarrow 0} \frac{1}{t}(\Phi(\mu_t)-\Phi(\mu)).$$
We have the following classical formula (see  \cite[Theorem 6.11]{Imayoshi:Taniguchi:book})
$$\dot\Phi_\mu([\nu])(z) = \left(-\frac{6}{\pi}\int\int_{f(\Half)}\frac{\lambda(\xi)}{(\xi-f(z))^4} |d\xi |^2\right) f'(z)^2 \qquad z \in \overline{\Half}$$
where $f = f_\mu$ and
$$\lambda(\xi) = \left( \frac{f_z}{\bar{f_z}}\frac{\nu}{1-|\mu|^2}\right)\circ f^{-1}(\xi).$$
We note that if $g_t = f_{\mu_t}\circ f_{\mu}^{-1}$ and $\lambda_t = \mu_{g_t}$ with
$$\lambda_t = t\lambda + O(t^2).$$
Further we note that $f_\mu$ conjugates $\Gamma$ to quasifuchsian group $\Gamma_\mu$ and if $\Omega = f_\mu(\Half)$ and $\Omega^* = f_\mu(\overline{\Half})$ then
$\Omega/\Gamma_\mu$ is conformal to $Y$. Furthermore $\lambda \in B(\Omega,\Gamma_\mu)$ is the beltrami differential on $Y$ corresponding to the deformation $\nu \in B(\Half,\Gamma)$.

Before we prove our bounds, we will prove a lemma which follows easily from the Area theorem.
\begin{lemma}
Let $\Omega,\Omega^*$ be complementary Jordan domains and $z \in \Omega$. Then 
$$\frac{1}{\rho_{\Omega}(z)} \int_{\Omega^*} \frac{|d\xi|^2}{|\xi-z|^4} \leq \frac{\pi}{4}.$$
where $\rho_\Omega$ is the hyperbolic metric on $\Omega$. 
\label{lemma:area}\end{lemma}

{\bf Proof:} If $M$ is a Mobius transformation then
$$|M(x)-M(y)|^2 = |M'(x)||M'(y)||x-y|^2$$
Thus if $M$ maps domains $\Omega_0,\Omega_0^*$ to $\Omega,\Omega^*$ and with $M(z_0) = z$  then
$$\frac{1}{\rho_{\Omega}(z)} \int_{\Omega^*} \frac{|d\xi|^2}{|\xi-z|^4} = \frac{1}{\rho_{\Omega_0}(z_0)} \int_{\Omega_0^*}\frac{|dw|^2}{|w-z_0|^4}.$$
Thus we can assume that $z=0$. Thus letting $\hat\Omega,\hat\Omega^*$ be the image of $\Omega, \Omega^*$ under $w= 1/\xi$. Then
$$\int_{\Omega^*} \frac{|d\xi|^2}{|\xi-z|^4}= \int_{\Omega^*} \frac{|d\xi|^2}{|\xi|^4} = \int_{\hat\Omega^*} |dw|^2= Area(\hat\Omega^*).$$
We now apply the area theorem. We choose $f:\Delta\rightarrow \Omega$ uniformizing $\Omega$ with $f(0) = 0$.
Then let $g(z) = 1/f(1/z)$. Then $g:\hat\Delta\rightarrow \hat\CC$ with complement of the image equal to $\hat\Omega^*$. By the area theorem (see \cite[Section II.1.5]{Lehto:book:univalent})
$$Area(\hat\Omega^*) \leq \frac{\pi}{|f'(0)|^2}.$$
Thus
$$\frac{1}{\rho_{\Omega}(z)} \int_{\Omega^*} \frac{1}{|\xi-z|^4}|d\xi|^2  \leq \frac{\pi}{|f'(0)|^2\rho_{\Omega}(z)} = \frac{\pi}{\rho_{\Delta}(0)} = \frac{\pi}{4}.$$
\eproof

Finally we have the following.
\begin{lemma}
The map $\Phi_X: \Teich(S) \rightarrow Q(\overline X)$ is $3/2$-Lipschitz with respect to the $L^\infty$ metric. Specifically if $u \in T_Y(\Teich(S))$
$$\|d\Phi_X(u)\|_\infty \leq \frac{3}{2}\|u\|_\infty.$$
\end{lemma}

{\bf Proof:}
$$\|\dot\Phi_\mu([\nu])\|_\infty =\sup_{z\in \overline{\Hp}} \frac{|\dot\Phi_\mu([\nu])(z)|}{\rho(z)}$$
Lifting we have
$$\|\dot\Phi_\mu([\nu])\|_\infty = \sup_{z\in \Omega^*} \frac{1}{\rho_{\Omega^*}(z)} \left| \frac{6}{\pi}\int_{\Omega} \lambda(\xi)\frac{1}{|\xi-z|^4}d\xi|^2\right|.$$
Thus
$$\|\dot\Phi_\mu([\nu])\|_\infty  \leq \frac{6\|\lambda\|_\infty}{\pi} \sup_{z\in \Omega^*} \left(\frac{1}{\rho_{\Omega^*}(z)} \int_{\Omega} \frac{1}{|\xi-z|^4}|d\xi|^2\right) .$$
By Lemma \ref{lemma:area} we then get
$$\|\dot\Phi_\mu([\nu])(z)\| \leq \frac{3}{2}\|\lambda\|_\infty.$$
\eproof

Although we do not need it, we also include the improved Lipschitz bound for the $L^2$ norm.

\begin{theorem}
The map $\beta_X: \Teich(S) \rightarrow Q(\overline X)$ is $3/2$-Lipschitz with respect to the $L^2$ metric. Specifically if $u \in T_Y(\Teich(S))$
$$\|d\beta_X(u)\|_2 \leq \frac{3}{2}\|u\|_2.$$
\end{theorem}

{\bf Proof:}  The tangent vector $u$ corresponds to path $\mu_t = \mu+t\nu+O(t^2)$ in $B(\Half,\Gamma) \simeq B(X)$ and path $\lambda_t = t\lambda+O(t^2)$ in $B(\Omega,\Gamma_\mu) \simeq B(Y)$.  Then  we have
$$\|\dot\Phi_\mu([\nu])\|_2^2 = \int_{\Half^*/\Gamma} \frac{|\dot\Phi_\mu([\nu])(z)|^2}{\rho(z)} |dz|^2 \leq \frac{6^2}{\pi^2}\int_{\Omega^*/\Gamma_\mu} \frac{1}{\rho_{\Omega^*}(z)} \left(\int_{\Omega} |\lambda(\xi)|\frac{1}{|\xi-z|^4}|d\xi|^2\right)^2 |dz|^2.$$
Applying Holder's inequality we have 
$$\|\dot\Phi_\mu([\nu])\|_2^2\leq \frac{36}{\pi^2}\int_{\Omega^*/\Gamma_\mu} \frac{|dz|^2}{\rho_{\Omega^*}(z)} \int_{\Omega} \frac{1}{|w-z|^4}|dw|^2 \int_{\Omega} |\lambda(\xi)|^2\frac{1}{|\xi-z|^4}|d\xi|^2.$$
By Lemma \ref{lemma:area},
 $$\|\dot\Phi_\mu([\nu])\|_2^2\leq \frac{9}{\pi}\int_{\Omega^*/\Gamma_\mu} |dz|^2 \int_{\Omega} |\lambda(\xi)|^2\frac{1}{|\xi-z|^4}|d\xi|^2.$$
On $\Omega\times\Omega^*$ we consider the area form
 $$\omega(\xi,z) =  \frac{|\lambda(\xi)|^2}{|\xi-z|^4}|dz|^2|d\xi|^2.$$
As Mobius transformations satisfy
$$| \gamma(w) - \gamma(z)|^2 = |\gamma'(z)||\gamma'(w)||z-w|^2$$  
then $\omega$ is invariant under the diagonal action of $\Gamma_\mu$ on $\Omega\times\Omega^*$. As both $(\Omega/\Gamma_\mu)\times\Omega^*$ and 
$\Omega\times(\Omega^*/\Gamma_\mu)$ are fundamental domains for the diagonal action we have
$$ \int_{\Omega\times(\Omega^*/\Gamma_\mu)} \omega = \int_{(\Omega\times\Omega^*)/\Gamma_\mu}\omega = \int_{(\Omega/\Gamma_\mu)\times\Omega^*} \omega.$$
It follows that 
$$\|\dot\Phi_\mu([\nu])\|_2^2\leq \frac{9}{\pi}\int_{\Omega^*} \frac{|dz|^2}{|\xi-z|^4} \int_{\Omega/\Gamma_\mu} |\lambda(\xi)|^2|d\xi|^2$$
By Lemma \ref{lemma:area}, we integrate again  to get
$$\|\dot\Phi_\mu([\nu])\|_2^2\leq \frac{9}{4} \int_{\Omega/\Gamma_\mu} |\lambda(\xi)|^2\rho_{\Omega}(\xi)|d\xi|^2  = \frac{9}{4}\|\lambda\|^2_2$$
\eproof

\bibliography{bib,math}
\bibliographystyle{math}

\end{document}